\documentclass[12pt]{amsart}
\usepackage{amsmath, amssymb}
\usepackage{enumerate} 
\makeatletter
\@namedef{subjclassname@2020}{\textup{2020} Mathematics Subject Classification}
\makeatother
\newtheorem{theorem}{\bf Theorem}[section]
\newtheorem{lemma}[theorem]{\bf Lemma}
\newtheorem{proposition}[theorem]{\bf Proposition}

\newtheorem{example}[theorem]{\sc Example}

\newcommand{\N}{\mathbb{N}}

\newcommand{\pr }{\mathrm{Pr} }

 \newcommand{\ol}{\overline}
\DeclareMathOperator{\aut}{Aut}
\DeclareMathOperator{\out}{Out}

\DeclareMathOperator{\Sym}{Sym}

\begin{document}
\title[Commuting probability]{Commuting probability for the Sylow subgroups of a profinite group}
\thanks{The first two authors are members of GNSAGA (INDAM), 
and the third author was  supported by  FAPDF and CNPq.
}

\author[E. Detomi]{Eloisa Detomi}
\address{Dipartimento di Matematica \lq\lq Tullio Levi-Civita\rq\rq, Universit\`a degli Studi di Padova, Via Trieste 63, 35121 Padova, Italy} 
\email{eloisa.detomi@unipd.it}
\author[M. Morigi]{Marta Morigi}
\address{Dipartimento di Matematica, Universit\`a di Bologna\\
Piazza di Porta San Donato 5 \\ 40126 Bologna \\ Italy}
\email{marta.morigi@unibo.it}
\author[P. Shumyatsky]{Pavel Shumyatsky}
\address{Department of Mathematics, University of Brasilia\\
Brasilia-DF \\ 70910-900 Brazil}
\email{pavel@unb.br}


\begin{abstract} 
Given two subgroups $H,K$ of a compact group $G$, the probability that a random element of $H$ commutes with a random element of $K$ is denoted by $\pr(H,K)$.

\noindent  We show that if $G$ is a profinite group containing a Sylow $2$-subgroup $P$, a Sylow $3$-subgroup $Q_3$ and a Sylow $5$-subgroup $Q_5$ such that $\pr(P,Q_3)$ and $\pr(P,Q_5)$ are both positive, then $G$ is virtually prosoluble (Theorem \ref{main1}). 

\noindent Furthermore, if $G$ is a prosoluble group in which for every subset $\pi\subseteq\pi(G)$ there is a Hall $\pi$-subgroup $H_\pi$  and a Hall $\pi'$-subgroup $H_{\pi'}$ such that $\pr(H_\pi,H_{\pi'})>0$, then  $G$ is virtually pronilpotent (Theorem \ref{main2}). 
 \end{abstract}
 \maketitle

\section{Introduction} 
If $H$ and $K$ are two subsets of a finite group, then the probability that two randomly chosen elements from $H$ and $K$ commute is given by
\[ \pr(H,K) =\frac{ | \{ (x,y) \in H\times K \mid xy=yx \} |}{|H|\,|K|}.\]

This concept can be naturally extended to compact groups, using the normalized Haar measure (see the next section for details). 

Information on the commuting probabilities of some specific subgroups  gives insights on the structure of the group, and recently there has been an increasing interest in this subject (see e.g. \cite{AS-compact, AS2, DS-commprob, DLMS-prob, DLMS-finite, DMS-approx, eberhard, gr, nath, compactprob}).

In \cite{DLMS-finite} the authors  studied the commuting probability of Sylow subgroups in finite groups. In particular, the following result was obtained. 
\\

\begin{emph}{
 Let $G$ be a finite group, and let  $\epsilon>0$ be a real number such that for any distinct primes  $p,q\in\pi(G)$  
 there exists   a Sylow $p$-subgroup $P$ and a Sylow $q$-subgroup $Q$ of $G$ such that $\pr(P,Q) \ge \epsilon$.
Then $G$ has a metanilpotent subgroup of $\epsilon$-bounded index. 
}
\end{emph}
\\

Throughout this paper, we say that a quantity is $t$-bounded if it is bounded by a number depending only on the parameter $t$. Recall that a group $K$ is metanilpotent if the quotient $K/F(K)$ over its Fitting subgroup is nilpotent. 

The above result can be naturally extended to profinite groups $G$ such that there is a positive number $\epsilon$ with the property that for any distinct primes  $p,q\in\pi(G)$ the group $G$ contains a Sylow $p$-subgroup $P$ and a Sylow $q$-subgroup $Q$ for which $\pr(P,Q) \ge \epsilon$ (see Proposition \ref{profinite-easy}). Here $\pi(G)$ denotes the set of primes dividing the order of the profinite group $G$, which is a supernatural number.

On the other hand, in the realm of profinite groups, it is more natural to handle the situation where the probabilities are positive, but not necessarily bounded away from zero. 

It is not difficult to see that a profinite group in which every pair of coprime Sylow subgroups has positive commuting probability is not necessarily virtually metapronilpotent (see e.g. Example \ref{no-meta}). Here and throughout, we say that a profinite group virtually has some property if it has an open subgroup with that property.
 
Nevertheless, in this paper we prove that the property of having positive commuting probabilities of Sylow subgroups has strong impact on the structure of a profinite group. Moreover, in some cases it is enough to look only at the commuting probabilities of Sylow subgroups for the primes $2, 3$ and $5$.

\begin{theorem}\label{main1}
 Let $G$ be a profinite group containing a Sylow $2$-subgroup $P$, a Sylow $3$-subgroup $Q_3$ and a Sylow $5$-subgroup $Q_5$ such that $\pr(P,Q_3)$ and $\pr(P,Q_5)$ are both positive. Then $G$ is virtually prosoluble. 
\end{theorem}

The proof of the above result uses the fact that a nonabelian finite simple group has even order which is divisible either by $3$ or by $5$. This is an immediate corollary of the Odd Order Theorem \cite{fetho} and Thompson's classification of minimal insoluble groups \cite{thompson}.

Amazingly, it seems Theorem \ref{main1} does not admit any quantitative analogue. In particular, the index of the prosoluble radical of $G$ is not bounded in terms of    $\pr(P,Q_3)$ and  $\pr(P,Q_5)$. Indeed, let $G_r= \textrm{PSL}(2,r)$, where $r\equiv 3,5 \textrm{ mod } 8$. The Sylow $2$-subgroup $P$ of $G_r$ is of order  $4$. Therefore $\pr(P,Q) > 1/4 $ for every Sylow subgroup $Q$ of $G_r$. On the other hand, the order of $G_r$ tends to infinity when $r$ does. 

Another result obtained in this paper provides a sufficient condition for a prosoluble group to be virtually pronilpotent. Namely, the following holds. 

\begin{theorem}\label{main2}
Let $G$ be a prosoluble group and in which for every subset $\pi\subseteq\pi(G)$ there is a Hall $\pi$-subgroup $H_\pi$  and a Hall $\pi'$-subgroup $H_{\pi'}$ such that $\pr(H_\pi,H_{\pi'})>0$. Then  $G$ is virtually pronilpotent. 
 \end{theorem} 

This should be compared with \cite[Theorem 1.2]{DLMS-finite} where it was shown that if $G$ is a finite soluble group for which there is $\epsilon>0$ such that $\pr(P,H_{p'})\ge\epsilon$ whenever $P$ is a Sylow $p$-subgroup and $H_{p'}$ is a Hall $p'$-subgroup, then $G$ has a nilpotent subgroup of $\epsilon$-bounded index.

The rest of the paper is divided into three parts. In the next section we develop machinery for the study of the commuting probability of subgroups in compact groups. Section \ref{tres} is essentially devoted to the proof of Theorem \ref{main1} and the last one to the proof of Theorem \ref{main2}.

\section{General facts on commuting probabilities} 

In this section we present some general results about commuting probabilities in compact groups. So, unless otherwise specified, $G$ will be a compact Hausdorff topological group and a subgroup of $G$ will always be a closed subgroup. Note that any finite index subgroup of $G$ is open in $G$. We recall  that a profinite group  is a  compact  Hausdorff topological space, where a basis of the neighbourhood of the identity consists of normal subgroups of finite index. 
 
The Borel $\sigma$-algebra $\mathcal M$ of a compact group $G$ is the one generated by all closed subsets of $G$. We say that a measure $\mu$ on $(G,\mathcal M)$ is a (left) Haar measure provided $\mu$ is both inner and outer regular, $\mu(K) < \infty$ and $\mu(xE) = \mu(E)$ for all compact subsets $K$ and measurable subsets $E$ of $G$ (see \cite[Chapter 4]{HR} or \cite[Chapter II]{Nach}). Recall that there is a unique Haar measure $\mu$ on $(G,\mathcal M)$ such that $\mu(G) = 1$. 

Moreover, if $H$ is a subgroup of $G$, then either $\mu(H) =0$ or $\mu(H) >0$, and in the latter case $H$ is open in $G$ and $\mu(H)=|G:H|^{-1}.$

If $H$ and $K$ are  subgroups of $G$ the
set 
\[C = \{(x, y) \in H \times K \mid xy = yx\}\]  is closed in $H \times K$ since it is the  preimage of $1$ under the continuous map $f : H\times K \rightarrow G$ given by $f(x, y) = [x, y]$.
Denoting the normalized Haar measures of $H$ and $K$ by $\mu_H$ and $\mu_K$, respectively, the probability that a random element from $H$ commutes with a random element from $K$ is defined as 
 \[\pr(H,K) = (\mu_H \times \mu_K)(C).\]

Note that 
\[\pr(H,K) = \int_H \mu_K (C_K(x)) d\mu_H(x)  = \int_K \mu_H (C_H(y)) d\mu_K(y).\]

For every $x \in G$, the centralizer $C_G(x)$ equals $f_x^{-1} (1)$, where $f_x$ is the continuous function $f_x(y) = [x, y]$. It follows that $C_G(x)$ is closed and measurable. 
We define 
\[ \pr (x, K)= \mu_K \left(C_K(x)\right), \]
and note that $\pr (x, K)= |K:C_K(x)|^{-1}$ whenever $C_K(x)$ is open in $K$. 
	
The following lemma is immediate from Remark 2.9 in  \cite{AS-compact}.
\begin{lemma}\label{Xclosed}
 Let $G$ be a compact group, $H, K$ subgroups of $G$ and $n$ a positive integer. 
 Then the  set $X = \{x \in H \mid |x^K| \le n \}$ is a closed subset of $H$.	
 \end{lemma}

\begin{lemma} 
Let $H,K$ be subgroups of a compact group $G$ and let 
\[X=\{ x\in H \mid \pr(x,K)>0 \}. \] 
Assume that $\mu_H(X)>0$. Then $\pr(H,K)>0$. 
\end{lemma}
\begin{proof}
For a positive integer $n$ set
	\[X_n:=\left\{x\in H \mid {\pr(x,K)} \geq {1}/{n}\right\}.\]
As $X_n=\left\{x\in H \mid |K:C_K(x)|\le n\right\},$ by  Lemma \ref{Xclosed}  it is a closed and therefore measurable subset of $H$. 
	  
Observe that $X$ is the union of the ascending chain of subsets
 $$X_2 \subseteq X_3 \subseteq X_4 \subseteq \dots.$$ Therefore 
	 \[\mu_H(X)=\mu_H \left( \bigcup_{n \in \mathbb N} X_n \right) = \sup_{n \in \mathbb N} \mu_H (X_n).\]
	  Since $\mu_H(X)>0$, there exists $n \in \N$ such that $\mu_H (X_n)>0$. 
	 We see that
	  \begin{align*}
	  \pr(H,K) & = \int_H \mu_K (C_K(x)) d\mu_H(x) 
	  \ge \int_{X_n} \mu_K (C_K(x)) d\mu_H(x) \\
	 &\ge  \int_{X_n} \frac{1}n d\mu_H(x) =\frac{\mu_H (X_n)}n>0.
	  \end{align*} 
	 This concludes the proof.
\end{proof}

The next lemma  is an obvious extension of Lemma 2.4 in \cite{AS-compact}. 
\begin{lemma}\label{sub} 
Let $H,K$ be subgroups of a compact group $G$. Then for any subgroup $H_0 \le H$ we have 
\[\pr(H_0,K)\geq \pr(H,K).\]
\end{lemma}
\begin{proof}
Given $x \in K$, the map 
\[\alpha: \{yC_{H_0}(x) \mid y \in H_0\} \rightarrow \{hC_{H}(x) \mid h \in H\}\] 
that takes $yC_{H_0}(x)$ to  $yC_{H}(x)$ 
 is injective. Hence, $|H_0: C_{H_0}(x)| \le |H:C_{H}(x)|$ and so 
 $\mu_{H_0}(C_{H_0}(x)) \ge \mu_{H}(C_{H}(x))$. 
 Therefore, 
\[\pr(H_0,K) =  \int_K \mu_{H_0} (C_{H_0}(x)) d\mu_K(x) \ge  \int_K \mu_{H} (C_{H}(x)) d\mu_K(x) = \pr(H,K),
\]
 as claimed. 
\end{proof}

The following lemma is useful for handling quotients and (normal) subgroups. It is a variant of Lemma 2.3 in \cite{DS-commprob} (see also  Lemma 2.4 in \cite{AS-compact}). 
\begin{lemma}\label{quot} 
 Let $G$ be a compact group and let $N$ be a normal subgroup of  $G$.  For any subgroups $H,K$ of $G$ we have 
\[\pr(H,K)\leq \pr(HN/N,KN/N)\pr(N\cap H,N\cap K).\]
\end{lemma}
\begin{proof} 
Let $x \in H$. If $\mu_K (C_K(x)) > 0$, then $C=C_K(x)$ is open in $K$ and $C(N\cap K)=CN \cap K$. Hence, 
\begin{align*}
 \mu_K (C)^{-1}&=|K:C|=|K:C(N\cap K)| \cdot |C(N\cap K):C|  \\
 &= |KN :CN| \cdot |N\cap K:N\cap C| \\
 &= \mu_{KN}(CN)^{-1}  \mu_{N\cap K}(N\cap C)^{-1}. 
\end{align*}
Thus, if $\Omega=\{ x \in H \mid \mu_K (C_K(x)) > 0\},$ then 
\begin{align*}
\pr(H,K) &= \int_H \mu_K (C_K(x)) d\mu_H(x) =  \int_\Omega \mu_K (C_K(x)) d\mu_H(x)  \\
 &= \int_\Omega \mu_{KN}(C_K(x) N) \cdot \mu_{N\cap K}(N\cap C_K(x)) d\mu_H(x)  \\
 &\le \int_H \mu_{KN}(C_K(x)N) \cdot \mu_{N\cap K}(C_{N\cap K}(x)) d\mu_H(x).
\end{align*}
By applying the extended Weil formula  (see \cite[p. 88]{reiter}), and setting $\ol K=KN/N$  we obtain
\begin{align*} 
& \pr(H,K)     \\
\nonumber &\le  \int_{\frac{H}{N\cap H}} \left( \int_{N\cap H}  \mu_{KN}(C_K(xh)N) \cdot \mu_{N\cap K}(C_{N\cap K}(xh)) d\mu_{N\cap H}(h) \right) d\mu_{\frac{H}{N\cap H}}(x(N\cap H))\\
\nonumber  &\le  \int_{\frac{H}{N\cap H}} \left( \int_{N\cap H}  \mu_{\ol K}(C_{\ol K}(xN)) \cdot \mu_{N\cap K}(C_{N\cap K}(xh)) d\mu_{N\cap H}(h) \right) d\mu_{\frac{H}{N\cap H}}(x(N\cap H))\\
\nonumber    &\le  \int_{\frac{H}{N\cap H}}  \mu_{\ol K}(C_{\ol K}(\ol x))  \left( \int_{N\cap H}  \mu_{N\cap K}(C_{N\cap K}(xh)) d\mu_{N\cap H}(h) \right) d\mu_{\frac{H}{N\cap H}}(x(N\cap H)).
\end{align*}
Now we want a bound for 
\[ \int_{N\cap H}  \mu_{N\cap K}(C_{N\cap K}(xh)) d\mu_{N\cap H}(h). \]
For each $x \in H$, define the set 
\begin{align*}
A_x &=\{(h,k) \in (N\cap H) \times (N\cap K) \mid [xh,k]=1\}\\
&=\{(h,k) \in (N\cap H) \times (N\cap K) \mid xh \in C_{H}(k) \cap x(N\cap H) \}.
\end{align*}

If $A_x  \neq \emptyset$, then there exists $y \in C_{H}(k) \cap x(N\cap H)$. We have $C_{H}(k) \cap y(N\cap H)= C_{ H}(k) \cap x(N\cap H)  $ and since  \[[yz,k]=[y,k]^z[z,k]=1 \ \text{if and only if} \ [z,k]=1,\]
we have that $C_{H}(k) \cap y(N\cap H)=yC_{N\cap H}(k)$. 
Thus 
\[ (h,k) \in A_x \text{ if and only if } h \in x^{-1} yC_{N\cap H}(k), \]
so $A_x= \{ (h,k) \in (N\cap H) \times (N\cap K) \mid h \in x^{-1} yC_{N\cap H}(k)\}$.
Now we use the Lebesgue-Fubini Theorem 
\begin{align*}
\int_{N\cap H}  \mu_{N\cap K}(C_{N\cap K}(xh)) d\mu_{N\cap H}(h) &=
\int_{(N\cap H)\times (N\cap K)} \chi_{A_x} (h,k)  d(\mu_{N\cap H} \times  \mu_{N\cap K})(h,k) \\
&\le \int_{N\cap K} \mu_{ N\cap H} (x^{-1} yC_{N\cap H}(k)   ) d \mu_{N\cap K}(k) \\
&=\int_{N\cap K} \mu_{ N\cap H} (C_{N\cap H}(k)   ) d \mu_{N\cap K}(k) \\
&= \pr(N\cap H, N\cap K).
\end{align*}
Replacing this in the main inequality we have 
\begin{align*}
\pr(H,K) \le  \pr(N\cap H, N\cap K) \cdot \int_{\frac{H}{N\cap H}}  \mu_{\ol K}(C_{\ol K}(\ol x))   d\mu_{\frac{H}{N\cap H}}(x(N\cap H)). 
\end{align*}
Finally, since $H/(N\cap H)$ is isomorphic to $HN/N$, it follows (see e.g. Corollary 2.5 in \cite{compactprob}) that
\begin{align*}
\pr(H,K) &\le  \pr(N\cap H, N\cap K) \cdot \int_{\frac{HN}{ H}}  \mu_{\ol K}(C_{\ol K}(\ol x))   d\mu_{\frac{HN}{H}}(xN) \\ 
 &= \pr(N\cap H, N\cap K) \cdot \pr(\ol H, \ol K),
\end{align*}
 as desired. 
\end{proof}

In what follows we write $\langle X \rangle$ to denote the subgroup generated by a set $X$.
The next lemma is essentially Lemma 2.1 in \cite{eberhard} (see \cite[Lemma 2.7]{AS-compact} for a  detailed proof). 
\begin{lemma}\label{<X>}
 Let $G$ be a compact group with normalized Haar measure $\mu$ and
let $r \ge 1$. Suppose that $X$ is a closed symmetric subset of $G$ containing the
identity. If $\mu (X) > (r+1)^{-1}$, then $\langle X \rangle = X^{3r}$.
 \end{lemma}

The following lemma will play a central role in the proofs of main theorems.

\begin{lemma}\label{lem:H_0}
Let $H, K$ be subgroups of a compact group $G$ for which $\pr (H,K) \ge \epsilon >0 $. Then there exists a subgroup $H_0$ of $H$ such that the following holds.  
 \begin{enumerate}
\item $|H:H_0| \le 2/\epsilon$; 
\item 
  $|K: C_K(x)| \le (2/\epsilon)^{6/\epsilon}$  for every $x \in H_0$.
\end{enumerate}
\end{lemma}
\begin{proof} 
Set
 \[X=\{x\in H \mid  |x^K|\leq 2/\epsilon\}.\]
 By Lemma \ref{Xclosed}, $X$ is closed, hence measurable. 
 Moreover 
 \[ H \setminus X= \{ h \in H \mid   |h^K|>  2/\epsilon \} =  \{ h \in H \mid   \mu_K(C_K(h)) <  \epsilon/2 \}. \]
 
 Therefore, 
 \begin{align*}
 \epsilon \le \pr (H,K) &= \int_H \mu_K(C_K(h)) d\mu_H(h) \\
  &= \int_X \mu_K(C_K(h)) d\mu_H(h) +  \int_{H \setminus X} \mu_K(C_K(h)) d\mu_H(h) \\
    &\le \int_X  d\mu_H(h) +  \int_{H \setminus X}  \frac{\epsilon}2 d\mu_H(h) \\
    &\le  \mu_H(X) +   \frac{\epsilon}2  (1- \mu_H(X)) \le  \mu_H(X) +   \frac{\epsilon}2  
 \end{align*}
 and so $\mu_H(X) \ge {\epsilon}/2$. 
 
   Let $H_0$ be the subgroup generated by $X$. 
Clearly,  $\mu_H (H_0 ) \ge  \mu_H (X) \ge {\epsilon}/2$, so $|H:H_0| \le 2/\epsilon$. 
 Moreover, it follows from Lemma  \ref{<X>} that every element of $H_0$ is a product of at most $6/\epsilon$ elements of $X$.
Therefore $|h^K| \le (2/\epsilon)^{6/\epsilon}$ for every $h \in H_0$.
\end{proof}

\section{Theorem \ref{main1}}\label{tres}  
We will start this section with some general comments on the commuting probability of Sylow subgroups in profinite groups.

We denote by $F_i(G)$ the $i$th term of the Fitting series of a profinite group $G$.  More specifically, $F_1(G)=F(G)$ is the Fitting sugroup of $G$, also called the pronilpotent radical of $G$. This is the product of the normal pronilpotent subgroups of $G$. For $i\ge1$ the subgroup $F_{i+1}(G)$ is defined as the inverse image of the Fitting subgroup of $G/F_i(G)$. Note that $G$ is metapronilpotent precisely when $G=F_2(G)$.

 As mentioned in the introduction,  Lemma \ref{quot} guarantees that Theorem 1.1 in \cite{DLMS-finite} can be extended to profinite groups. Indeed, let $G$ be a  profinite group,  $\epsilon$ a positive real number, and assume that for any distinct primes $p$, $q\in\pi(G)$ there is a Sylow
$p$-subgroup $P$ and a Sylow $q$-subgroup $Q$ of $G$ such that $\pr(P, Q) \ge\epsilon$. 
 By Lemma \ref{quot}, for any open normal subgroup $N$ of $G$ we have 
 \[\pr(PN/N,QN/N) \ge \epsilon. \]
Hence, by  \cite[Theorem 1.1]{DLMS-finite}, the quotient group $G/N$ has a metanilpotent normal subgroup of index at most $n$, where $n$ is a positive integer depending on $\epsilon$ only.  Since $G$ is an inverse limit of the groups $G/N$, we conclude that $G$ contains a metapronilpotent subgroup of index at most $n$. 

This proves the following proposition. 

\begin{proposition}\label{profinite-easy} Let $\epsilon$ be a positive real number and $G$ a  profinite group in which for any distinct primes $p$, $q\in\pi(G)$ there is a Sylow $p$-subgroup $P$ and a Sylow $q$-subgroup $Q$ such that $\pr(P, Q) \ge\epsilon$. Then $G$ has a metapronilpotent subgroup of $\epsilon$-bounded index. 
\end{proposition}

\smallskip

One would be tempted to guess
 that if the profinite group $G$ satisfies the weaker hypothesis that for any distinct primes $p$, $q\in\pi(G)$ there is a Sylow
$p$-subgroup $P$ and a Sylow $q$-subgroup $Q$ of $G$ such that $\pr(P, Q)>0,$ then $G$ is virtually 
 metapronilpotent. However the following example shows that this is false.

\begin{example}\label{no-meta} {\rm Let $\{p_1, p_2, \dots\}$ be an infinite set of primes, and let $P_i$ be the cyclic group of order $p_i$. We assume that  $p_i<p_{i+1}$  for every $i$. Consider the iterated wreath product $G_i=P_i\wr\dots\wr P_2\wr P_1$
 and set $G$ to be the inverse limit of the groups $G_i$.  Note that, for every $i$, a $p_i$-Sylow subgroup $S_{p_i}$ of $G$ is finite. Furthermore,
  $\pr(S_{p_i},S_{p_j})\ge 1/|S_{p_i}|>0$ but $G$ is not virtually metapronilpotent.}
\end{example}
\medskip
 
In the rest of this section, we deal with the proof of Theorem \ref{main1}.
We will use the well-known fact that if $G$ is profinite group,  $H$ is a closed subgroup of $G$ and  $H_0$ is an open subgroup of $H$ then there exists an open normal subgroup $M$ of $G$ such that $M\cap H\le H_0.$ 

The following result is Theorem 2.2 in \cite{BFMMNSST}.
\begin{theorem}\label{TAMS} Let $X$ be a finite nonabelian simple group of order divisible by an
odd prime $r$. Then $X$ contains a conjugacy class of $r$-elements of even size.
\end{theorem}

We also need a technical lemma.

\begin{lemma}\label{direct-product}  Let $K=S_1\times\dots\times S_t$ be the direct product of finite simple groups ($t$ factors). Assume that $r$ is an odd prime dividing the order of $S_i$ for every $i$.  Let $P$ be a Sylow $2$-subgroup and $Q$ a Sylow $r$-subgroup of $K$ such that 
 \[ |P:C_P(x) | \le \eta \ \textrm{ for every }  x \in Q\]
for some positive number $\eta.$ Then $t$ is $\eta$-bounded.
\end{lemma}
\begin{proof}  Write $P=P_1\times\dots\times P_t$ and $Q=Q_1\times\dots \times Q_t$, where $P_i$ (resp. $Q_i$) is a Sylow $2$-subgroup (resp. $r$-subgroup) of $S_i$ for each $i$. By Theorem \ref{TAMS}, for each $i$ there exists $x_i\in Q_i$ such that the conjugacy class of $x_i$ has even size. Therefore $C_{P_i}(x_i)$ has index at least $2$ in $P_i$. Taking $x=(x_1,\dots,x_t)$ it follows that $|P:C_P(x)|\ge 2^t$. Therefore $2^t\le \eta$ and this implies the result.
\end {proof} 

It is well known that every finite group $K$ possesses a series
\[1 = K_0 \le K_1 \le\dots \le K_{2h+1} = K\]
of normal subgroups such that $K_{i+1}/K_i$ is soluble (possibly trivial) if $i$ is
even and a direct product of nonabelian simple groups if $i$ is odd. Following \cite{KhSh15} the number
of insoluble factors in this series is called the insoluble length $\lambda(K)$ of $K$.  Wilson showed in \cite{Wilson} that if $\mathcal{X}_1,\dots, \mathcal{X}_n$ are classes of finite groups closed with  respect to normal subgroups and subdirect products and 
if $\mathcal{X}$ is the class of finite groups having a normal series of given length $n$ such that
the $i$-th section is a $\mathcal{X}_i$-group, then any pro-$\mathcal{X}$ group has a normal series of length $n$  such that
the $i$-th section is a pro-$\mathcal{X}_i$-group. In particular, a combination of Lemma 2 and Lemma 3 of
\cite{Wilson} shows that if $\mathcal{X}$ is the class of finite groups $K$ such that $\lambda(K) \le l,$ then
any pro-$\mathcal{X}$ group has a normal series of  length at most $2l+1$ each of whose factors is
either prosoluble or a Cartesian product of nonabelian finite simple groups.

Recall that the generalized Fitting subgroup $F^*(G)$ of a finite group
$G$ is the product of the Fitting subgroup $F(G)$ and all subnormal quasisimple subgroups. Here a group is quasisimple if it is perfect and its
quotient by the centre is a nonabelian simple group.

The proof of the next lemma uses the  Schreier conjecture, i.e. the fact that  the outer automorphism groups of finite simple groups are soluble, which is a consequence of the classification of finite simple groups. We  also recall here that all nonabelian finite simple $3'$-groups are Suzuki groups, and the order of any Suzuki group is divisible by $5$.

\begin{lemma}\label{lambda} Let $\eta$ be a positive number and assume that  $G$ is a finite group  containing a Sylow $2$-subgroup $P$, a Sylow $3$-subgroup $Q_3$ and a Sylow $5$-subgroup $Q_5$ such that 
 \[ |P:C_P(x) | \le \eta \ \textrm{ for every }  x \in Q_3 \cup Q_5 .\]
Then $\lambda(G)$ is $\eta$-bounded.
\end{lemma} 
\begin{proof}  Without loss of generality, we can assume that the soluble radical of $G$ is trivial. Then $F^{*}(G)$ is a direct product of finite simple groups. 
As the order of a finite simple group is either divisible by $3$ or by $5$, we can assume that $F^{*}(G)=S_1\times\dots\times S_k$, where 
$S_1,\dots,S_t$ are simple groups of order divisible by $3$ and $S_{t+1},\dots,S_k$ of order divisible by $5$. In view of Lemma \ref{direct-product} applied to the group $S_1\times\dots\times S_t$ with $r=3$ it follows that $t$ is $\eta$-bounded.  Similarly, observe that $k-t$ is $\eta$-bounded. Hence, $k$ is $\eta$-bounded.

Note that $C_G(F^*(G)) \le Z(F^*(G))$ (see, for instance, \cite[Theorem 9.8]{Isaacs}). Since $Z(F^*(G))=1$, it follows that $G$ is isomorphic to a subgroup of the automorphism group $\aut(F^*(G))$, which in turn can be embedded in the semidirect product $B\rtimes \Sym(k)$, where the symmetric group $\Sym(k)$ acts on $B=\aut(S_1)\times\dots\times\aut(S_k)$ by permuting the factors (here the action is not necessarily transitive). 

Therefore $\bar G=G/F^*(G)$ is isomorphic to a subgroup of $\bar B\rtimes\Sym(k),$ where $\bar B=\out(S_1)\times\dots\times\out(S_k)$. For each $i$, the outer automorphism group $\out(S_i)$ is soluble by the Schreier conjecture. Hence, $\lambda(\bar G)\le \lambda(\bar B\bar G/\bar B)$. Since $\bar B\bar G/\bar B$ has order at most $k!$, it follows that $\lambda(G)$ is $\eta$-bounded, as desired.
\end{proof} 

We are now ready to prove Theorem \ref{main1}. We will make use of the well-known fact that if a profinite group $G$ is finite-by-prosoluble, then it is virtually prosoluble. Indeed, if $K$ is a finite normal subgroup of $G$ such that $G/K$ is prosoluble, then $C_G(K)$ is an open prosoluble subgroup.

\begin{proof}[Proof of Theorem \ref{main1}]
Let $\epsilon>0$ be a real number and $G$ a profinite group containing a Sylow $2$-subgroup $P$, a Sylow $3$-subgroup $Q_3$ and a Sylow $5$-subgroup $Q_5$ such that $\pr(P,Q_3) \ge \epsilon$ and $\pr(P,Q_5) \ge \epsilon$. We wish to prove that $G$ is virtually prosoluble. 
  
 By Lemma \ref{lem:H_0}, there exist two subgroups $R_3 \le Q_3$ and $R_5 \le Q_5$ such that, for  $i=3,5$,  $|Q_i :R_i| \le 2/\epsilon$ and $|P: C_P(x)| \le (2/\epsilon)^{6/\epsilon}$ for every $x \in R_i$. 
  
 As $R_i$ is open in the induced topology of $Q_i$, there exists an open normal subgroup $M_i$ of $G$ such that $Q_i\cap M_i\le R_i$. The subgroup $M=M_3\cap M_5$ is open and normal in $G$. Observe that $M \cap Q_i \le R_i$. Therefore for $i=3,5$ the intersection $M\cap Q_i$ is a Sylow $i$-subgroup of $M$ with the property that  $|P: C_P(x)| \le (2/\epsilon)^{6/\epsilon}$ for every $x \in M \cap Q_i$. It is clear that $P \cap M$ is a  Sylow $2$-subgroup of $M$ and $$|P \cap M: C_{P \cap M}(x)| \le |P: C_P(x)|.$$ 
	
Therefore  it is sufficient to prove that $M$ is virtually prosoluble. We can simply assume that $G=M$ and
 \[ (*)\ \ \ \ \ \ \ \ |P:C_P(x) | \le \eta \ \textrm{ for every }  x \in Q_3 \cup Q_5, \] where $\eta$ is a number depending on $\epsilon$ only.
 
 Note that the condition $(*)$ is inherited by homomorphic images and normal subgroups of $G$. 
 By Lemma \ref{lambda} the insoluble length $\lambda(G/N)$ of every finite continuous image $G/N$ of $G$ is $\eta$-bounded. 
 As already mentioned, 
 it follows from \cite{Wilson} that  $G$ has a normal series of  $\eta$-bounded length each of whose factors is
either prosoluble or a Cartesian product of nonabelian finite simple groups. 

Let $T$ be a factor in this series of the form  $T=\prod_{i\in I} S_i$, where each $S_i$ is a nonabelian finite simple group. For every finite subset $J$ of $I$ the subgroup $\prod_{i\in J} S_i$ is normal in $T$. It follows from  Lemma \ref{direct-product} that  $|J|$ is  $\eta$-bounded. Therefore the cardinality of  $I$ is $\eta$-bounded and in particular $T$ is finite. 

Now the conclusion follows by induction on $\lambda(G)$. Namely, if $\lambda(G)=0$ then $G$ is prosoluble and the result holds. Assume that 
$\lambda(G)=h\ge 1$ and let $R$ be the maximal normal prosoluble subgroup of $G$ (the prosoluble radical). Passing to the quotient $G/R$, without loss of generality we can assume that $R=1$. In view of the above $G$ possesses a finite normal subgroup $T$, which is a direct product of nonabelian simple groups, such that $\lambda(G/T)=h-1$. By induction, $G/T$ is virtually prosoluble. It follows that $G$ is virtually prosoluble as well and  the proof is complete.
\end{proof}

\section{Theorem \ref{main2}} 

This section is devoted to the proof of Theorem \ref{main2}. We deal with profinite groups $G$ satisfying the following condition
\begin{itemize}
\item[(**)] For every subset $\pi\subset\pi(G)$ there exists a Hall $\pi$-subgroup $H_\pi$ and a Hall $\pi'$-subgroup $H_{\pi'}$ such that $\pr(H_\pi,H_{\pi'})>0$. 
\end{itemize}

Remark that in view of Lemma \ref{quot} the condition $(**)$ is inherited by quotients and normal subgroups of $G$.

Our first goal will be to prove that a nontrivial profinite group $G$ satisfying $(**)$ has nontrivial Fitting subgroup. The following elementary lemma is needed. 

\begin{lemma}\label{fc}
 If $G$ is a prosoluble group containing a nontrivial element $x\in G$ such that $|G:C_G(x)|$ is finite,  then $F(G)\ne 1$.
\end{lemma} 
\begin{proof}
Assume by contradiction that $F(G)=1$.  Observe that $K=\langle x^G\rangle$ is generated by finitely many conjugates of $x$ and so  its centralizer has finite index in $G$. It follows that the centre $Z(K)$ of $K$ has finite index in $K$. Moreover $Z(K)\le F(G)=1$. Therefore $Z(K)=1$ and consequently $K$ is finite. As $G$ is prosoluble, $K$ is soluble and $F(K)\ne 1$. Since $F(K)\le F(G)$,  this contradicts the assumption that $F(G)=1$.
\end{proof}

\begin{lemma}\label{pq} 
 Let $\epsilon$ a positive number and $G$ a profinite group containing subgroups $L,K$ such that $\pr (L,K) \ge \epsilon$. Let $p,q$ be primes bigger than $(2/\epsilon)^{6/\epsilon}$. Then every $p$-element of $L$ commutes with every $q$-element of $K$. 
\end{lemma} 
\begin{proof} 
Let $P$ be an arbitrary Sylow $p$-subgroup of $L$ and $Q$ be an arbitrary Sylow $q$-subgroup of $K$.  By Lemma \ref{sub}, $\pr(P,Q) \ge \pr (L,K)  \ge \epsilon$.  Lemma \ref{lem:H_0} shows that there is an open subgroup  $Q_0\leq Q$ such that  $|Q:Q_0| \le 2/\epsilon$ and  $|P: C_P(y)| \le (2/\epsilon)^{6/\epsilon}$  for every $y \in Q_0$. 
     Since $p, q > (2/\epsilon)^{6/\epsilon}$, we have that $Q=Q_0$ and also  $P= C_P(y)$ for every $y\in Q$. Thus $[P,Q]=1$. 
\end{proof}

\begin{proposition}\label{Fitting}
Let $G$ be a nontrivial prosoluble group and assume that  for every subset $\pi\subset\pi(G)$ there exists a Hall $\pi$-subgroup $H_\pi$ 
 and a Hall $\pi'$-subgroup $H_{\pi'}$ such that $\pr(H_\pi,H_{\pi'})>0$.  
 Then  $F(G)\neq1$. 
 \end{proposition} 
 
\begin{proof}
Suppose that $G$ has an infinite Sylow $p$-subgroup $P$. Let $H_{p'}$ be a Hall $p'$-subgroup of $G$ such that $\pr(P,H_{p'})>0$. 
By Lemma \ref{lem:H_0}, $P$ has an open subgroup $P_0\leq P$ with the property that, whenever $x\in P_0,$ the centralizer $C_{H_{p'}}(x)$ has finite index in $H_{p'}$. Note that $P_0\ne 1$ since $P$ is infinite.

Choose a nontrivial element $x \in P_0$.  As  $C_{H_{p'}}(x)$ is open in $H_{p'}$,  there is an open normal subgroup $N$ in $G$ such that $N\cap  H_{p'} \le C_{H_{p'}}(x)$. Write $N_{p'}=N \cap  H_{p'}$ and $N_p=N \cap P$. Note that $N_{p'}$ is a  $p'$-Hall subgroup  of $N$ contained in  $C_{N}(x)$.

We have $N=N_{p'}N_p$. For an arbitrary element $h\in N$ write $h=h_{p'}h_p $, where  $h_{p'} \in N_{p'}$ and $h_p \in N_p$. Thus, $x^h=x^{h_p}\in P$ for every $h\in N$. It follows that $[N,x]\le P$  is a pro-$p$ group. Note that $[N,x]$ is normal in $N$. So if $[N,x] \neq 1$, then $F(N)\neq1$, whence $F(G)\neq1$. 
 
Otherwise, $N$ is contained in the centralizer of $x$ in $G$. It follows that $|G:C_G(x)|$ is finite, and $F(G)\ne 1$ by Lemma \ref{fc}. 
 
So we can assume that all Sylow subgroups in $G$ are finite. If there exists an element $x$ in $G$ such that $\pi(|G:C_G(x)|)$ is finite, then the index $|G:C_G(x)|$ is finite and we deduce from  Lemma \ref{fc} that $F(G)\ne1$.  

Therefore we  are left with the case where  $\pi(|G:C_G(x)|)$ is infinite for any $x\in G$ and we will show that this assumption leads to a contradiction. 
Note that in particular $\pi(G)$ is infinite. 

We will choose two disjoint infinite subsets $\pi=\{p_1,p_2,\dots\}$ and $\sigma=\{q_1,q_2,\dots\}$ of $\pi(G)$ such that for every $i$ there is a $p_i$-element $x_i$ such that $q_i\in\pi(|G:C_G(x_i)|)$. This can be done inductively in the following manner.

Let $p_1$ be any prime divisor of the order of $G$. Let $x_1$ be any nontrivial $p_1$-element, and let $p_1\neq q_1\in\pi(|G:C_G(x_1)|)$. Set $\pi_1=\{p_1\}$ and $\sigma_1=\{q_1\}$. 
 
 Further, assume we have chosen the sets $\pi_i=\{p_1,\dots,p_i\}$ and  $\sigma_i=\{q_1, \dots , q_i\}$ with the required properties. Taking into account that $\pi(G)$ is infinite pick  
 $p_{i+1} \in \pi(G)\setminus \{ p_1, \dots , p_i,   q_1, \dots , q_i\}$ and a nontrivial $p_{i+1} $-element $x_{i+1} $. Since $\pi(|G:C_G(x_{i+1})|)$ is infinite, we can choose  a prime $q_{i+1} \in \pi(|G:C_G(x_{i+1})|)\setminus \{ p_{i+1}, p_1, \dots , p_i,   q_1, \dots , q_i\}$.
 
 Observe that the sets $\pi=\cup\pi_i$ and $\sigma=\cup\sigma_i$ are as required.

 By hypothesis, $G$ contains a Hall $\pi$-subgroup $K$  and a Hall $\pi'$-subgroup $L$ such that $\pr(K,L)=\epsilon$ for some positive number $\epsilon$.
 By Lemma \ref{pq},  whenever $p\in\pi(K)$,  $q\in\pi(L)$ and $p,q\geq (2/\epsilon)^{6/\epsilon}$, 
  any $p$-element of $K$ commutes with any $q$-element  of $L$. Since both $\pi$ and $\sigma$ are infinite, there exists $i$ such that both $p_i,q_i\ge (2/\epsilon)^{6/\epsilon}$. By replacing $x_i$ with a suitable conjugate, we can assume that  $x_i\in K$. This yields a contradiction because  $\sigma\subseteq\pi(L)$ meaning that $L$ contains a Sylow $q_i$-subgroup of $G$ and so $q_i\not\in\pi(|G:C_G(x_i)|).$
\end{proof}

The following result is well known. It easily follows from  \cite[Lemma 4.29]{Isaacs} using the standard inverse limit argument. 

\begin{lemma} \label{coprime} Let $G$ be a profinite group, $N$ a normal subgroup of $G$ and $A$ a subgroup of $G$ such that $\pi(N)$ and $\pi(A)$ are disjoint. Then $[N, A, A]=[N, A].$\end{lemma}

We will also need the following technical result, which is Lemma 2.6 in \cite{DMScoprime}. 

\begin{lemma}\label{2.6 coprime} Assume that $G = QH$ is a finite group with a normal nilpotent subgroup $Q$ and a subgroup H such that $(|Q|, |H|) = 1$ and $|Q:C_Q(x)| \le m$ for all $x\in H.$ Then, the order of $[Q, H]$ is $m$-bounded.
\end{lemma}

We are now ready to prove Theorem \ref{main2}. Recall that the pronilpotent residual $M$ of a profinite group $G$ is the intersection of all normal open subgroups $N$ of $G$ such that $G/N$ is nilpotent.

\begin{proof}[Proof of Theorem \ref{main2}]
Let $G\ne 1$ be a prosoluble group in which for every subset $\pi\subset\pi(G)$ there is a Hall $\pi$-subgroup $H_\pi$ 
 and a Hall $\pi'$-subgroup $H_{\pi'}$ such that $\pr(H_\pi,H_{\pi'})>0$.  
We need to prove that $G$ is virtually pronilpotent.
 
It follows from Proposition \ref{Fitting} that $F(G)\neq 1$. Note that, by Lemma \ref{sub}, $F_2(G)$ satisfies the hypothesis of theorem. 

Suppose first that $G=F_2(G)$.
 
 Let $M$ be the pronilpotent residual of $G$, choose $p\in\pi(M)$, and let $P$ be a Sylow $p$-subgroup of $M$. 
 Note that $P$ is normal in $G$ since $M$ is pronilpotent. 
 By hypothesis there exists  a Sylow $p$-subgroup $G_p$ and  a Hall $p'$-subgroup $H_{p'}$  of $G$ such that $\pr(G_p, H_{p'})=\epsilon$ 
  for some positive number $\epsilon$. 
  Since $P \le G_p$, by Lemma \ref{sub} we have that $\pr( P, H_{p'}) \ge\epsilon$. 
 Further, by Lemma \ref{lem:H_0} there exists a subgroup  $H_0$  of $H_{p'}$ such that $| H_{p'}:H_0| \le 2/\epsilon$ and
   $| P: C_{ P}(x)| \le (2/\epsilon)^{6/\epsilon}$  for every $x \in H_0$. 

We now wish to show that the following holds.
\begin{center}$(*)$$\qquad\qquad H_{p'}/C_{H_{p'}}(P)$ is finite.\end{center}

Set $R=PH_0$. By Lemma \ref{2.6 coprime} the order of $[P,H_0]$ is $\epsilon$-bounded in every finite continuous image of $R$. Therefore $[P,H_0]$ is finite and subsequently $C=C_{H_0}([P,H_0])$ is open in $H_0$. By Lemma \ref{coprime},   $[P, C]=[P, C,C]\le [P,H_0, C]=1$
so $C \le C_{H_0}(P)$. 
It follows that $C_{H_0}(P)$ has finite index in $H_{p'}$. This proves $( *)$. 

In the case where $\pi(M)$ is finite we prove that $G$ is virtually pronilpotent using induction on $|\pi(M)|$.

 If $|\pi(M)|=1$, then  $M$ is a pro-$p$ subgroup. 
Let $G_p$  be  Sylow $p$-subgroup of   $G$. Clearly $G_p$ contains $M$.  Since $G/M$ is pronilpotent, $H_{p'}$ is pronilpotent and  $G_p$  is normal in $G$. 
 Let  $V=G_pC_{H_{p'}} (M)$. 
  As $H_{p'}/C_{H_{p'}}(M)$ is finite by $(*)$, the subgroup $V$ is open in $G$. Moreover, $[G_p, C_{H_{p'}}(M)] \le M$. We deduce that $[G_p, C_{H_{p'}}(M),C_{H_{p'}}(M)] =1$. In view of Lemma \ref{coprime} this implies that $[G_p, C_{H_{p'}}(M)] =1$, whence $V$ is pronilpotent, and $G$ is virtually pronilpotent. 
  
So assume that $|\pi(M)|\geq2$. Write $M=M_1\times M_2$, where $M_1$ and $M_2$ are of coprime orders. By induction, both $G/M_1$ and $G/M_2$ are virtually pronilpotent. It follows that $G$ is virtually pronilpotent as well, being isomorphic to a subgroup of $G/M_1\times G/M_2$.

From now on we will therefore assume that $\pi(M)$ is infinite. For every $p \in \pi(M)$, let $S_p$ be a Sylow $p$-subgroup of $M$ and set $\sigma_p=\pi(H_{p'}/C_{H_{p'}}(S_p))$.    Note that in view of $(*)$  each  $\sigma_p$ is finite. Moreover, $\sigma_p\neq\emptyset$ since $S_p=[S_p,H_{p'}]$ (see for example \cite[Lemma 2.4]{AST}).

Suppose that the union $\Sigma=\cup_{p \in \pi(M)}\sigma_p$ is infinite. In this case we can choose two disjoint infinite sets of primes $\pi_0=\{p_1,p_2,\dots\}\subseteq\pi(M)$ and $\sigma_0=\{q_1,q_2,\dots\}\subseteq\pi(G)$ such that $q_i\in\sigma_{p_i}$  for every index $i$. To see this, we argue as follows. First, choose any $p_1\in\pi(M)$ and $q_1\in\sigma_{p_1}$. Then assume we have found primes $p_1, \dots , p_i$ and $q_1, \dots , q_i$ with the required properties. Since $\Sigma$ is infinite,  there is a prime $q_{i+1}\in \Sigma $ such that $q_{i+1}\not\in\sigma_{p_1}\cup \dots\cup \sigma_{p_i}\cup\{ p_1, \dots , p_i\}$. Then $q_{i+1}\in\sigma_t$ for some prime $t \in \pi(M)\setminus \{ p_1, \dots , p_i\} $, and we set $p_{i+1}=t$. In this way we obtain the required sets $\pi_0$ and $\sigma_0$. 

By hypothesis, $G$ contains a Hall $\pi_0$-subgroup $K$ and a Hall $\pi'_0$-subgroup $L$ such that $\pr(K,L)=\epsilon$ for some positive number $\epsilon$. 
  It follows from  Lemma \ref{pq} that for any $p\in\pi(K)$, $q\in\pi(L)$ such that $p,q\geq (2/\epsilon)^{6/\epsilon}$ the $p$-elements of $K$ commute with the $q$-elements of $L$. Since both $\pi_0$ and $\sigma_0$ are infinite, there exists $i$ such that $p_i,q_i\ge (2/\epsilon)^{6/\epsilon}$. 
This yields a contradiction since, 
 on the one hand,  the $p_i$-Sylow subgroup $S_{p_i}$ of $M$ is contained in $K$ and $q_i \in \sigma_0\subseteq\pi(L)$ while, on the other hand, $q_i\in\pi(|G:C_G(S_{p_i})|)$.

So we will assume that the union $\Sigma$   is finite. Let $N=C_G(M)$. Note that $N/(M\cap N)$ is pronilpotent, being isomorphic to a subgroup of $G/M$. Moreover, $M\cap N$ is contained in the center of $N$ and so $N$ is pronilpotent. As $G/N$ is isomorphic to a subgroup of the direct product $\prod_{p \in \pi(M)}(G/C_G(S_p))$, it follows that $\pi(G/N)\subseteq \cup_p\sigma_p=\Sigma$ is finite. In particular,  $\pi(G/F(G))$ is finite since $N\le F(G)$.
	
For every $q \in \pi(G/F(G))$, there exists a Sylow $q$-subgroup $Q$ and a Hall $q'$-subgroup $H_{q'}$ of $G$ such that $\pr (Q, H_{q'}) >0$. We wish to show that $Q \cap F(G)$ has finite index in $Q$. 
		    
Without loss of generality we may assume that  $O_q(G)=1$, in which case
 $F(G)$ is a $q'$-subgroup contained in any Hall $q'$-subgroup of $G$.
   It follows from Lemma \ref{sub}  that $\pr(Q,F(G))>0$.
 By Lemma \ref{lem:H_0}  there exists an open subgroup $Q_0$ of $Q$ 
 such that  $|F(G): C_{F(G)}(x)| \le (2/\epsilon)^{6/\epsilon}$  for every $x \in Q_0$. 
 By  Lemma \ref{2.6 coprime},   the order of $[F(G),Q_0]$ is $\epsilon$-bounded in every finite continuous image of $G$ and so we deduce that $[F(G),Q_0]$ is finite. Thus, $C=C_{Q_0}([F(G),Q_0])$ has  finite index in $Q$. Moreover 
  \[[F(G),C,C]\le[F(G),Q_0,C]=1.\]
By Lemma \ref{coprime}, this implies that $C$ centralizes $F(G)$.  Recall that the Fitting subgroup of a metapronilpotent group contains its centralizer. Since  $F(G)$ is a $q'$-group, it follows that $C_Q(F(G))=1$ and so $Q$ is finite. 

  Thus, we have shown that for any $q\in\pi(G/F(G))$ and any Sylow $q$-subgroup $Q$ the index of $Q\cap F(G)$ in $Q$ is finite. 
  Since $\pi(G/F(G))$ is finite, we deduce that $G/F(G)$ is finite. This proves that  if $G$ is metapronilpotent, then $G$ is virtually pronilpotent. 
  
 Now we drop the assumption that $G$ is metapronilpotent.  By the  above, $F_2(G)/F(G)$ is finite. So $G$ contains an open normal subgroup $U$ such that  $U\cap F_2(G) \leq F(G)$. Since $F_2(G)/F(G)=F(G/F(G))$, this implies that $UF(G)/F(G)$ has trivial Fitting  subgroup. By Lemma \ref{quot} the quotient group $UF(G)/F(G)$ satisfies the hypothesis of  Proposition \ref{Fitting}. We conclude that $UF(G)/F(G)=1$.  Therefore $U\leq F(G)$ has finite index in $G$ and the group $G$ is virtually pronilpotent. 
\end{proof}

\end{document}